\documentclass[12pt]{amsart}
\usepackage{latexsym}
\usepackage{amstext}
\usepackage{amssymb}
\usepackage{amsfonts}
\usepackage{amsthm}
\usepackage{amsopn}
\usepackage{amsbsy}
\usepackage{layout}
\usepackage{color}
\usepackage{amsmath}
\usepackage{graphicx}
\usepackage{bm}
\usepackage{yfonts}

\newtheorem{deff}{Definition}[section]
\newtheorem{lemma}[deff]{Lemma}
\newtheorem{theorem}[deff]{Theorem}
\newtheorem{corollary}[deff]{Corollary}

\newtheorem{proposition}[deff]{Proposition}

\newtheorem{em-example}[deff]{Example}
\newtheorem{em-def}[deff]{Definition}        
\newtheorem{em-remark}[deff]{Remark}         
\newtheorem{em-question}[deff]{Question}

\newtheorem{problem}[deff]{Problem}

\newenvironment{example}{\begin{em-example} \em }{ \end{em-example}}

\newenvironment{remark}{\begin{em-remark} \em }{\end{em-remark}}

\newcommand{\acal}{{\mathcal A}}

\newcommand{\dcal}{\Delta}

\newcommand{\fcal}{\mathcal {F}}

\newcommand{\pcal}{\mathcal {P}}

\def\ker{\mathop{\rm ker}}

\def\cl{\mathop{\it cl}}

\DeclareMathSymbol{\res}{\mathord}{AMSa}{"16}

\def\:{\nobreak \hskip .1111em\mathpunct {}\nonscript \mkern
   -\thinmuskip {:}\hskip .3333emplus.0555em\relax}

\catcode`\@=12

\def\N{{\mathbb N}}
\def\R{{\mathbb R}}
\def\Q{{\mathbb Q}}

\begin{document}
\title[Basic properties of $X \in \dcal$]{Basic properties of $X$ for which spaces $C_p(X)$ are distinguished}

\author{Jerzy K\c akol}
\address{ Faculty of Mathematics and Informatics, A. Mickiewicz University,
61-614 Pozna\'{n}, Poland and Institute of Mathematics Czech Academy of Sciences, Prague, Czech Republic}
\email{kakol@amu.edu.pl}
\author{Arkady Leiderman}
\address{Department of Mathematics, Ben-Gurion University of the Negev, Beer Sheva, Israel}
\email{arkady@math.bgu.ac.il}

\begin{abstract}
In our paper \cite{KL} we showed that a  Tychonoff space $X$ is a $\Delta$-space (in the sense of \cite{Knight}, \cite{Reed}) if
and only if the locally convex space $C_{p}(X)$ is distinguished.
Continuing this research, we investigate whether the class $\Delta$ of $\Delta$-spaces is invariant under the basic topological operations.

We prove that if $X \in \Delta$ and $\varphi:X \to Y$ is a continuous surjection such that $\varphi(F)$ is an $F_{\sigma}$-set in $Y$ for every closed set $F \subset X$,
then also $Y\in \Delta$. As a consequence, if $X$ is a countable union of closed subspaces $X_i$ such that each $X_i\in \Delta$, then also $X\in \Delta$.
 In particular, $\sigma$-product of any family of scattered Eberlein compact spaces is a $\Delta$-space and
the product of a $\Delta$-space with a countable space is a $\Delta$-space.
Our results give answers to several open problems posed in \cite{KL}.

Let $T:C_p(X) \longrightarrow C_p(Y)$ be a continuous linear surjection. We observe that $T$ admits an extension to a linear continuous operator $\widehat{T}$
from $\R^X$ onto $\R^Y$ and deduce that $Y$ is a $\Delta$-space whenever $X$ is.
 Similarly, assuming that $X$ and $Y$ are metrizable spaces, we show that $Y$ is a $Q$-set whenever $X$ is.

Making use of obtained results,  we provide a very short proof for the claim that every compact $\Delta$-space has countable tightness.
As a consequence, under  Proper Forcing Axiom (PFA) every compact $\Delta$-space is sequential.

In the article we pose a dozen open questions.
\end{abstract}

\thanks{The research for the first named author is supported  by the GA\v{C}R project 20-22230L and RVO: 67985840.
He thanks  also the Center For Advanced Studies in Mathematics of Ben-Gurion University of the Negev for financial support during his visit in 2019.}

\maketitle

\section{Introduction}\label{intro}
Throughout the article, all topological spaces are assumed to be Tychonoff.
By $C_{p}(X)$ we mean the space of all real-valued continuous functions on a Tychonoff space $X$
 endowed with the topology of pointwise convergence.

\begin{deff}{\rm (\cite{KL}, \cite{Knight})}\label{def:Delta}
A topological space $X$ is said to be a $\Delta$-space if for every decreasing sequence $\{D_n: n \in \omega\}$
of subsets of $X$ with empty intersection, there is a decreasing sequence $\{V_n: n \in \omega\}$ consisting of open
subsets of $X$, also with empty intersection, and such that $D_n \subset V_n$ for every $n \in \omega$.
\end{deff}

The class of all $\Delta$-spaces is denoted by $\Delta$.
Let us point out that the original definition of a $\Delta$-set $X \subset \R$, where $\R$ denotes the real line,
is due to G. M. Reed and E. K. van Douwen (see \cite{Reed}).
$\Delta$-sets of reals have been used and investigated thoroughly in the study of two of the most basic and central constructions in general topology:
the Moore--Nemytskii plane and the Pixley-Roy topology.
Denote by $M(X)$ the subspace of the Moore--Nemytskii plane, which is obtained by using only a subset $X \subset \R$ of the $x$-axis.
G. M. Reed observed that $M(X)$ is countably paracompact if and only if $X$ is a $\Delta$-set \cite{Reed}.

For a separable metrizable space $M$, denote by $\fcal(M)$ the hyperspace of finite subsets of $M$ endowed with the Pixley-Roy topology.
D. J. Lutzer proved that if $M$ is a strong $\Delta$-set, i.e. every finite power $M^n$ is a $\Delta$-set, then $\fcal(M)$ is countably paracompact \cite{Lutzer}.
H. Tanaka proved the converse statement: if $\fcal(M)$ is countably paracompact, then $M$ is a strong $\Delta$-set \cite{Tanaka}.
Also, the work \cite{Tanaka} deals with the analogous questions for general (not necessarily separable) metrizable spaces.

A set of reals $X$ is called a $Q$-set if each subset of $X$ is $F_{\sigma}$, or, equivalently,
each subset of $X$ is $G_{\delta}$ in $X$. The existence of uncountable $Q$-sets is independent of ZFC. Every $Q$-set is a $\Delta$-set, but consistently the converse is not true
(see\cite{Knight}). More details about $Q$-sets and $\Delta$-sets can be found in \cite{FM}, \cite{Knight}.
Of course, there are plenty of nonmetrizable $\Delta$-spaces with non-$G_{\delta}$ subsets, in ZFC \cite{KL}.

We could not find a single paper devoted to investigation of the general topological $\Delta$-spaces. 
Quite recently the authors have shown that the notion of $\Delta$-spaces plays a key role in the study of distinguished $C_p$-spaces \cite{KL}.
\begin{theorem}\cite{KL}\label{Theor:description}
For a Tychonoff space $X$, the following conditions are equivalent{\rm:}
\begin{enumerate}
\item[{\rm (1)}] $C_p(X)$ is distinguished.
\item[{\rm (2)}] Any countable disjoint collection of subsets of $X$ admits a point-finite open expansion in $X$.
\item[{\rm (3)}] $X$ is a $\Delta$-space.
\end{enumerate}
\end{theorem}
We should mention that independently and simultaneously an analogous description of distinguished $C_p$-spaces (but formulated in different terms)
appeared in \cite{FS}.

By a \emph{bounded set} in a topological vector space  we understand any set which is absorbed by every $0$-neighbourhood.
Following  J. Dieudonn\'{e} and
L. Schwartz  \cite{dieudonne} a locally convex space (lcs) $E$ is called \emph{distinguished} if every bounded subset of the bidual of $E$
 in the weak$^{*}$-topology is contained in the closure of the weak$^{*}$-topology  of  some bounded subset of $E$. 
Equivalently, a lcs $E$ is distinguished if and only if the strong dual of $E$  (i.e. the topological  dual of $E$ endowed with the strong topology) 
is \emph{barrelled}. A. Grothendieck  \cite{grothendieck} proved that a metrizable lcs $E$ is distinguished 
if and only if its strong dual is \emph{bornological}. 
Recall that the \emph{strong topology} on $E^{\prime}$ is the topology of uniform convergence on bounded subsets of $E$.

Denote by $L_p(X)$ the dual of $C_p(X)$, i.e. the linear space of all continuous linear functionals on $C_p(X)$, endowed with the topology of pointwise convergence.
Basic properties of $L_p(X)$ are described thoroughly in \cite{Arch5}.
By $L_{\beta}(X)$ we denote the strong dual of $C_{p}(X)$, i.e. the space $L_p(X)$ endowed with the \emph{ strong topology} $\beta_{X}=\beta(L_p(X), C_{p}(X))$.
 Note also that for a vector space $E$ the finest locally convex topology $\xi$ of $E$ is generated by the family of  all absolutely convex and absorbing subsets of $E$ 
which form a base of neighbourhoods of zero for the topology $\xi$. 

The following main characterization theorem has been proved instrumental in the study of distinguished lcs $C_p(X)$.\\
\\

\begin{theorem}[\cite{FKLS}, \cite{fe-ka}]\label{three}
For a Tychonoff space $X$ the following assertions are equivalent:
\begin{enumerate}
\item[{\rm (1)}] $C_{p}(X)$ is distinguished, i.e. its strong dual $L_{\beta}(X)$ is a barrelled space.
\item[{\rm (2)}] For each $f\in\mathbb{R}^{X}$ there is a bounded  $B\subset C_{p}(X)$ with $f\in \cl_{\R^X}(B)$.
\item[{\rm (3)}] The strong topology $\beta_{X}$ of the strong dual of $C_{p}(X)$ is the finest locally convex topology on $L_{p}(X)$.
\end{enumerate}
\end{theorem}

Naturally, aforementioned crucial Theorem \ref{Theor:description} has been proved in \cite{KL} with the help of Theorem \ref{three}.
In this paper, Theorem \ref{three} has been applied effectively again for the proof of Theorem \ref{linear1}, the main result of Section \ref{Linear}.

Our aim is to continue the research about topological $\Delta$-spaces originated in our paper \cite{KL}.
We obtain results in two directions.
First, in Section \ref{Images} we investigate whether the class $\Delta$ is invariant under the basic topological operations, including
continuous images, closed continuous images, countable unions and finite products.
\emph{What do we know about continuous images of $\Delta$-spaces?}

\begin{proposition}\cite{KL} \label{prop:locomp}
There exists in ZFC a MAD family $\acal$ on $\N$ such that the corresponding Isbell--Mr\'owka
space $\Psi(\acal)$ admits a continuous mapping onto the closed interval $[0,1]$.
\end{proposition}
Thus, the class $\Delta$ is not invariant under continuous images even for first-countable separable locally compact pseudocompact spaces.
The following result has been proved in our paper \cite{KL}.
\begin{proposition}\cite{KL} \label{prop:closedmap} Let $X$ be any $\Delta$-space and $\varphi: X \to Y$ be a closed continuous surjection with finite fibers.
Then $Y$ is also a $\Delta$-space.
\end{proposition}

Shortly after the paper \cite{KL} was published, V. Tkachuk \cite{Tkachuk} observed that the proof of Proposition \ref{prop:closedmap} in fact
does not use the last restriction about finiteness of fibers.
 So, Proposition \ref{prop:closedmap} is valid without unnecessary assumption of finiteness of fibers and the class $\Delta$ is invariant under closed continuous images.
As an immediate consequence, V. Tkachuk \cite{Tkachuk} noticed that we have a positive answer to Problem 5.3 posed in \cite{KL}:
any continuous image of a compact $\Delta$-space is also a $\Delta$-space.

In this paper we generalize Proposition \ref{prop:closedmap} as follows:
 Let $X$ be any $\Delta$-space and $\varphi: X \to Y$ be a continuous surjection
such that $\varphi(F)$ is an $F_{\sigma}$-set in $Y$ for every closed set $F \subset X$;
then $Y$ is also a $\Delta$-space (Theorem \ref{th:F-sigma_map}).
 It is interesting to note that the proof of Theorem \ref{th:F-sigma_map}
is obtained by absolutely elementary arguments.

We say that a topological space $X$ is \emph{$\sigma$-closed discrete} if $X = \bigcup_{n\in\omega}X_n$, where each $X_n$ is a closed and discrete subset of $X$.
It is easy to see that every $\sigma$-closed discrete space is in $\Delta$.
A straightforward application of Theorem \ref{th:F-sigma_map} gives a far-reaching generalization of this fact:
Assume that $X$ is a countable union of closed subsets $X_n$, where each $X_n\in \Delta$; then also $X\in \Delta$ (Proposition \ref{prop:count_unions}).
As a corollary we solve in the affirmative Problem 5.8 posed in \cite{KL}:
a countable union of compact $\Delta$-spaces is also a $\Delta$-space.
In particular, $\sigma$-product of any family consisting of scattered Eberlein compact spaces is a $\Delta$-space.
Another consequence says that the product of a $\Delta$-space with a $\sigma$-closed discrete space
(in particular, a countable space) is a $\Delta$-space.
Remark that the general question whether the class $\Delta$ is invariant under finite products remains open.
It is worthwhile mentioning that we do need an assumption on finite fibers for the following "reverse" version of Proposition \ref{prop:closedmap} above:
Let $\varphi: X \to Y$ be a continuous surjection with finite fibers; then $Y \in \Delta$ implies that also $X \in \Delta$ (Proposition \ref{prop:reverse}).

Following A. V. Arkhangel'skii \cite{Arch4}, we say that a space $Y$ is $l$-dominated ($u$-dominated, $t$-dominated)
by a space $X$ if $C_p(X)$ can be mapped linearly and continuously (uniformly continuously, continuously, respectively) onto $C_p(Y)$.
There are many topological properties which are invariant under defined above relations, and there are many which are not.
The main goal of Section \ref{Linear} is to study the following question: \emph{Which topological properties related to being a $\Delta$-space are preserved by
the relation of $l$-dominance}?

We show that the class of Tychonoff spaces $\Delta$ is invariant under the relation of $l$-dominance, equivalently, the class of distinguished $C_p$-spaces 
is invariant under the operation of taking continuous linear images (Theorem \ref{linear1}).
For the reader's benefit, aiming to emphasize a big potential in this research area,
we present two different proofs of Theorem \ref{linear1}. The first proof is based on item (2) of Theorem \ref{three} and invokes a new observation
about extensions of linear continuous surjections between $C_p$-spaces. The second proof uses the language of dualities and is based on item (3) of Theorem \ref{three}.
As an easy consequence of our approach, we obtain that inside the class of subsets of $\R$, $Q$-sets are invariant under the relation of $l$-dominance (Corollary \ref{Q-set}).
What makes the proof of Corollary \ref{Q-set} surprising is the fact that it does not depend on the question whether the square of a $Q$-set remains a $Q$-set.
We prove also that such topological properties as $\sigma$-scattered, $\sigma$-discrete, Eberlein compact, scattered Eberlein compact are preserved by the relation of $l$-dominance.  

The last Section \ref{PFA} is devoted to the study of compact / countably compact $\Delta$-spaces.
It has been shown in \cite{LS} that every compact $\Delta$-space has countable tightness.
Relying on Proposition \ref{prop:closedmap}, we  provide a very short proof of this assertion.
As a consequence, under  Proper Forcing Axiom (PFA) every compact $\Delta$-space is sequential.
In \cite{KL} we proved that every compact $\Delta$-space $X$ is \emph{scattered}, i.e. every subset $A$ of $X$ has an isolated (in $A$) point.
Making use of Proposition \ref{prop:closedmap} again, in Theorem \ref{th:count_compact} we generalize this result for countably compact $\Delta$-spaces, in ZFC.

In order to better understand the boundaries of the class $\Delta$ one should search for new examples and counter-examples.
Another recent relevant paper \cite{LS} is devoted to the following question: under what conditions 1) trees topologies;
2) $\Psi$-spaces built on maximal almost disjoint families of countable sets; 3) ladder system spaces do belong to the class $\Delta$.
Note that there are compact scattered spaces $X \notin \Delta$ (for example, the compact space $[0,\omega_{1}]$) \cite{KL}.
A stronger result has been obtained in \cite{LS}: there exists a compact scattered space $X$ such that the scattered height of $X$ is finite, and yet $X \notin \Delta$.
Thus, Problem 5.11 from \cite{KL} has been solved negatively in \cite{LS}.

Our notations are standard, the reader is advised to consult with the monographs \cite{Arch5} and \cite{Engelking} for the notions which are not explicitly defined in the text.
In the article we pose a dozen open questions.
\section{Continuous images, unions and products of $\Delta$-spaces}\label{Images}

\begin{theorem}\label{th:F-sigma_map} Let $X$ be a $\Delta$-space and $\varphi: X \to Y$ be a continuous surjection
such that $\varphi(F)$ is an $F_{\sigma}$-set in $Y$ for every closed set $F \subset X$.
Then $Y$ is also a $\Delta$-space.
\end{theorem}
\begin{proof}
Assume that $\{D_n: n \in \omega\}$ is any decreasing sequence of subsets of $Y$ with empty intersection.
 Then $\{\varphi^{-1}(D_n): n \in \omega\}$ is a decreasing sequence of subsets of $X$ and
$\bigcap_{n\in\omega}\varphi^{-1}(D_n) = \varphi^{-1}(\bigcap_{n\in\omega}D_n) = \emptyset$.
By assumption, there is a decreasing sequence of open sets $\{U_n: n \in \omega\}$ such that $\varphi^{-1}(D_n) \subset U_n$ for each $n \in \omega$ and $\bigcap_{n\in\omega} U_n = \emptyset$.
Define $H_n = Y\setminus \varphi(X\setminus U_n)$ for each $n\in\omega$.
We have that $D_n \subset H_n$ for each $n\in\omega$. Indeed, $\varphi^{-1}(y) \subset U_n$ for every $y \in D_n$, hence
$y \notin \varphi(X\setminus U_n)$ which means that $y \in H_n$.
Note that $\varphi^{-1}(H_n) \subset U_n$, therefore $\bigcap_{n\in\omega} H_n = \emptyset$
 because $$\varphi^{-1}(\bigcap_{n\in\omega} H_n) = \bigcap_{n\in\omega} \varphi^{-1}(H_n) \subset  \bigcap_{n\in\omega} U_n = \emptyset$$
Clearly, the sets  $X\setminus U_n$ are closed, the sets $\varphi(X\setminus U_n)$ are $F_{\sigma}$ in $Y$, hence the sets $H_n$ are $G_{\delta}$ in $Y$.

It appears that we can refine the sets $H_n$ further and construct another decreasing sequence consisting of open sets $\{V_n: n \in \omega\}$ such that $H_n \subset V_n$ for each $n \in \omega$ and
 $\bigcap_{n\in\omega} V_n = \emptyset$. This claim in fact should be attributed to E. K. van Douwen (see \cite[page 150]{Reed}).
Because the real argument has not been provided in \cite{Reed}, for the sake of completeness we include the proof. Denote  by $H_n^k$ open subsets of $Y$ such that $H_{n+1}^k \subset H_n^k$
and $\bigcap_{n\in\omega}H_n^k = H_k$ for every $n, k \in \omega$. Since the sequence $\{H_k: k \in \omega\}$ is decreasing, by induction over upper index $k$ without loss of generality
we may assume also that $H_{n}^{k+1} \subset H_n^k$ for every $n, k \in \omega$. We declare now that $V_n = H_n^n$, $n \in \omega$. It is clear that $D_n \subset V_n$ and each $V_n$ is an open set.
We show that the intersection of all sets $V_n$ is empty. Indeed, let $y$ be any element of $Y$. There exists $k\in\omega$ such that $y \notin H_k$, hence there exists $n(k) \in \omega$ such that
$y \notin H_{n(k)}^k$. Fix any $m \geq \max\{k, n(k)\}$. Then $H_m^m \subset H_m^k \subset H_{n(k)}^k$ which implies that $y \notin H_m^m$. Finally,
$\bigcap_{n\in\omega}V_n = \emptyset$ and the proof is complete.
\end{proof}

\begin{proposition}\label{prop:count_unions}
Assume that $X$ is a countable union of closed subsets $X_n$, where each $X_n$ belongs to the class $\Delta$.
Then $X$ also belongs to $\Delta$. In particular, a countable union of compact $\Delta$-spaces is also a $\Delta$-space.
\end{proposition}
\begin{proof} Denote by $Z$ the free topological union of the spaces $X_n, n \in \omega$. It is easy to see that $Z \in \Delta$, by Theorem \ref{Theor:description}.
The space $Z$ admits a natural continuous mapping $\varphi$ onto $X$.
Since $\varphi(F)$ is an $F_{\sigma}$-set in $X$ for every closed set $F \subset Z$, we deduce that $X \in \Delta$, by Theorem \ref{th:F-sigma_map}.
 \end{proof}

Thus, we have a positive solution of Problem 5.8 posed in \cite{KL}.

\begin{corollary} \label{cor:cont_comp} Let $X$ be a $\sigma$-compact $\Delta$-space and $Y$ be a continuous image of $X$.
Then $Y$ also is a $\Delta$-space.
\end{corollary}

Thus, we have a positive solution of Problem 5.3 posed in \cite{KL}.

\begin{corollary}\label{cor:sigma_product}
$\sigma$-product of any family consisting of scattered Eberlein compact spaces is a $\Delta$-space.
\end{corollary}
\begin{proof} $\sigma$-product is a countable union of $\sigma_n$-products, where $\sigma_n$-product includes elements of the product
which support consists of at most $n$ points, $n \in \omega$ . Every $\sigma_n$-product of scattered Eberlein compact spaces is again a scattered Eberlein compact,
therefore it is a $\Delta$-space \cite{FKLS}. It remains to apply Proposition \ref{prop:count_unions}.
\end{proof}

\begin{corollary}\label{Lin} Let $X$ be a Lindel\" of subspace of a $\sigma$-product of any family consisting of scattered Eberlein compact spaces
and $Y$ be a continuous image of $X$. Then $Y$ also is a $\Delta$-space.
\end{corollary}
\begin{proof} $X$ is equal to the countable union of its closed subspaces $X_n$, where $X_n$ is the intersection of $X$ with $\sigma_n$-product.
Every Lindel\" of subspace of a scattered Eberlein compact is necessarily $\sigma$-compact, by a recent result of V. Tkachuk \cite{Tkachuk2}.
Finally, $Y$ is a $\dcal$-space by Corollary \ref{cor:cont_comp}.
\end{proof}

\begin{corollary}\label{Lin_Cech} Let $X$ be a Lindel\" of \v Cech-complete $\Delta$-space
and $Y$ be a continuous image of $X$. Then $Y$ also is a $\Delta$-space.
\end{corollary}
\begin{proof} Any Cech-complete $\Delta$-space is scattered \cite{KL}. Now we use the well-known fact stating that every Lindel\" of \v Cech-complete scattered
space is $\sigma$-compact (see \cite[Theorem 4.5]{Aviles}) and we finish the proof again by Corollary \ref{cor:cont_comp}.
\end{proof}

We don't know  answers to the following problems.

\begin{problem}\label{prob1}
Let $X$ be any Lindel\" of subspace of a compact $\Delta$-space. Is $X$ a $\sigma$-compact space?
\end{problem}

In the case that the answer to Problem \ref{prob1} is negative we can ask

\begin{problem}\label{prob2}
Let $X$ be any Lindel\" of subspace of a compact $\Delta$-space and $Y$ be a continuous image of $X$.
Is $Y$ a $\Delta$-space?
\end{problem}

We have a partial positive result for products of $\Delta$-spaces.

\begin{corollary}\label{product}
Let $Z$ be the product of a $\Delta$-space $X$ with with a $\sigma$-closed discrete space
(in particular, a countable space) $Y$. Then $Z$ also is a $\Delta$-space.
\end{corollary}
\begin{proof} Let $Y= \bigcup_{n\in\omega} Y_n$, where each $Y_n$ is a closed and discrete subset of $Y$.
Denote by $Z_n = X \times Y_n$. It is clear that each $Z_n$ is closed in $Z$ and $Z_n \in \Delta$.
We get that $Z$ is a countable union of closed $\Delta$-spaces $Z_n$, so Proposition \ref{prop:count_unions} applies.
\end{proof}

Next statement formally is more general than Proposition \ref{prop:count_unions}.

\begin{proposition}\label{prop:sigma-loc_finite}
Assume that $X$ is covered by a $\sigma$-locally finite family of closed subsets $\{X_{\alpha}: \alpha \in A\}$, where every $X_{\alpha}$ belongs to $\Delta$.
Then $X$ also belongs to $\Delta$.
\end{proposition}
\begin{proof} The union of a locally finite family of closed subsets of $X$ is closed in $X$ \cite{Engelking}.
Now remind the following fact \cite[Theorem 2.7]{Burke}. Suppose that $\pcal$ is a topological property preserved
under closed mappings and $\{X_{\alpha}: \alpha \in A\}$
is a locally finite closed cover of $X$ with each $X_{\alpha}$ satisfying $\pcal$. If the
free topological sum satisfies $\pcal$ then so does $X$. It suffices to say that $X$ satisfies property $\pcal$ if $X\in \Delta$ and apply Proposition \ref{prop:count_unions}.
\end{proof}

\begin{remark} Proposition \ref{prop:count_unions} is not valid without assuming that all pieces $X_n$ in the union are closed.
 Let $M$ be the Michael line which is the refinement of the real line $\R$ obtained by isolating all irrational points.
Clearly, $M$ can be represented as a countable disjoint union of singletons (rationals) and an open discrete set.
 Nevertheless, the Michael line $M$ is not in $\Delta$ \cite{FKLS}.
\end{remark}

Now we consider a question of "reversing" of Proposition \ref{prop:closedmap}.
It is evident that if $\varphi: X \to Y$ is a continuous one-to-one mapping from $X$ onto $Y$ and $Y$ is a $\Delta$-space, then $X$ is also a $\Delta$-space.
The following more general result has been conjectured by V. Tkachuk \cite{Tkachuk} and below we provide a straightforward argument.

\begin{proposition} \label{prop:reverse} Let $\varphi: X \to Y$ be a continuous finite-to-one surjective mapping. 
If $Y$ is a $\Delta$-space, then $X$ is also a $\Delta$-space.
\end{proposition}
\begin{proof}
Let $\{B_n: n \in \omega\}$ be any decreasing sequence of subsets of $X$ with empty intersection.
Denote by $D_n = \varphi(B_n)$. Then $\{D_n: n \in \omega\}$ is a decreasing sequence of subsets of $Y$ with empty intersection.
Indeed, for every $y\in Y$ the fiber $\varphi^{-1}(y)$ is finite, hence there is $n \in \omega$
such that $\varphi^{-1}(y) \cap B_n = \emptyset$ which means that $y \notin D_n$.
By assumption, there is a decreasing sequence of open sets in $Y$,  $\{U_n: n \in \omega\}$ such that $D_n \subset U_n$ for each $n \in \omega$ and $\bigcap_{n\in\omega} U_n = \emptyset$.
Define $V_n = \varphi^{-1}(U_n)$ for each $n\in\omega$. Clearly, a decreasing sequence of open in $X$ sets, $\{U_n: n \in \omega\}$ is as required.
\end{proof}

We don't know under which conditions the latter Proposition \ref{prop:reverse} can be generalized for the mappings with countable fibers.

\section{$\Delta$-spaces vs properties of spaces $C_{p}(X)$}\label{Linear}
Our main goal here is to study the following question: \emph{Which topological properties related to being a $\Delta$-space are preserved by 
the relation of $l$-dominance}?

The class of all distinguished lcs does not preserve continuous linear images.
To see this it suffices to consider the identical mapping from the Banach space $C[0,1]$ onto $C_p[0,1]$.
Below we show that the class of Tychonoff spaces $\Delta$ is invariant under the relation of $l$-dominance, equivalently, the class of distinguished $C_p$-spaces 
is invariant under the operation of taking continuous linear images. 

\begin{theorem}\label{linear1}
Assume that $Y$ is $l$-dominated by $X$.
 If $X$ is a $\Delta$-space, then $Y$ also is a $\Delta$-space.
\end{theorem}
For the reader's benefit, we present two different proofs of Theorem \ref{linear1}: topological and analytical ones.
In order to present the first proof, we start with the following simple lemma. Surprisingly, we were unable to find its formulation in any monograph
cited in the references. 
By this reason we include its complete proof which rely on several extreme properties of the Tychonoff product $\R^X$.

\begin{lemma}\label{help}
Let $X$ and $Y$ be two sets and let $E \subset \R^X$ and $F\subset \R^Y$ be dense vector subspaces of $\R^X$ and $\R^Y$, respectively.
Assume that $T: E \longrightarrow F$ is a continuous linear surjection between lcs $E$ and $F$.
 Then $T$ admits a continuous linear surjective (unique) extension $\widehat{T}:\R^X \longrightarrow \R^Y$. 
\end{lemma}
\begin{proof} Let us list all well-known properties of $\R^X$ we are going to use.\newline
\underline{Property 1}. Every closed vector subspace $H$ of $\R^X$ is complemented in $\R^X$ and the quotient
$\R^X/H$ is linearly homeomorphic to the product $\R^Z$ for some set $Z$ \cite[Corollary 2.6.5, Theorem 2.6.4]{bonet}.\newline
\underline{Property 2}. The product topology on $\R^X$ is minimal, i.e. $\R^X$ does not admit a weaker Hausdorff locally convex topology \cite[Corollary 2.6.5(i)]{bonet}. \newline
\underline{Property 3}. $\R^Y$ fulfills the extension property, i.e. if $M$ is a vector subspace of a lcs $L$, then every continuous linear mapping $T: M \longrightarrow \R^Y$
admits a continuous linear extension $\widehat{T}: L \longrightarrow \R^Y$ \cite[Theorem 10.1.2 (a)]{Narici}.
 
By Property 3, there exists a continuous linear extension $\widehat{T}: \R^X \longrightarrow \R^Y$ of $T$ such that $F\subset \widehat{T}(\mathbb{R}^{X})$.
 We prove that $\widehat{T}$ is a surjective mapping.
Denote by $\varphi:\mathbb{R}^{X}/\ker(T)\longrightarrow\mathbb{R}^{Y}$ the injective mapping associated with  the quotient mapping $Q: \mathbb{R}^{X}\longrightarrow\mathbb{R}^{X}/\ker(\widehat{T})$,
 where $\ker(\widehat{T})$ is the kernel of $\widehat{T}$ and $\varphi\circ Q = \widehat{T}$.
By Property 1, the space $\mathbb{R}^{X}/\ker(T)$ is linearly homeomorphic to the product $\mathbb{R}^{Z}$ for some set $Z$. 
So we may assume that $\varphi$ is  a continuous linear bijection from $\mathbb{R}^{Z}$ onto a dense subspace $\widehat{T}(\mathbb{R}^{X})$ of $\mathbb{R}^{Y}$.
  This implies that on $\widehat{T}(\mathbb{R}^{X})$ there exists a stronger locally convex topology $\xi$ such that $(\widehat{T}(\mathbb{R}^{X}),\xi)$
	is linearly homeomorphic with $\mathbb{R}^{Z}$.
	However, by Property 2, $\mathbb{R}^{Z}$ does not admit  a weaker Hausdorff locally convex topology, hence $\widehat{T}(\mathbb{R}^{X})$ 
	is isomorphic to the complete lcs $\mathbb{R}^{Z}$. Finally, $\widehat{T}(\mathbb{R}^{X})$ is closed in $\mathbb{R}^{Y}$ and then $\widehat{T}$ is a surjection.
 \end{proof}
\begin{proof}[First Proof]
Let $T:C_{p}(X)\longrightarrow C_{p}(Y)$ be a continuous linear surjection. Denote by
 $\widehat{T}:\mathbb{R}^{X}\longrightarrow \mathbb{R}^{Y}$ the extension of $T$ which is supplied by Lemma \ref{help}.
By Theorem \ref{Theor:description}, $C_{p}(X)$ is distinguished and we can apply item (2) of Theorem \ref{three}.
Take arbitrary $f\in\mathbb{R}^{Y}$. There exists $g\in\mathbb{R}^{X}$ with $\widehat{T}(g)=f$. Then there exists a bounded set $B\subset C_{p}(X)$ 
such that  $g\in \cl_{\R^X}(B)$. We define $A = T(B)$. It is easy to see that $A$ is bounded and 
$f\in \cl_{\R^Y}(A)$ which means that $C_{p}(Y)$ is distinguished, equivalently, $Y$ is a $\Delta$-space, by Theorem  \ref{Theor:description}. 
\end{proof}

\begin{proof}[Second Proof]
If $T:C_{p}(X)\longrightarrow C_{p}(Y)$ is a continuous linear surjection, then by \cite[Proposition 23.30, Lemma 23.31]{meise},
the adjoint mapping  $T^*: (L_p(Y),\beta_{Y})\longrightarrow (L_p(X),\beta_{X})$ is continuous and injective,
 where $\beta_{X}$ and $\beta_{Y}$ are the strong topologies on the duals $L_{p}(X)$ and $L_{p}(Y)$, respectively. 
Denote by $Z=T^*(L_p(Y))$. Endow $Z$ with the induced topology $\beta_X\hskip-2.5pt\restriction_Z$. Since
$T^*: (L_p(Y),\beta_{Y})\rightarrow (Z,\beta_{X}|Z)$ is a continuous linear bijection, the sets $T^*(U)$, where $U$ run over all absolutely convex neighbourhoods of zero in 
$(L_p(Y),\beta_{Y})$, form a base of absolutely convex neighbourhoods of zero  for a locally  convex topology $\xi$ on $X$ such that  $\beta_{X}\hskip-2.5pt\restriction_Z\leq \xi$ and
 $T^*: (L_p(Y),\beta_{Y})\longrightarrow (Z,\xi)$ is a  linear homeomorphism.
Since $C_{p}(X)$ is distinguished by Theorem \ref{Theor:description}, the topology  $\beta_X$  is the  finest locally convex topology, by item (3) of Theorem \ref{three}.
 The property of having the finest locally convex topology is inherited by vector subspaces, so the induced topology  $\beta_{X}\hskip-2.5pt\restriction_Z$  is the finest locally convex one. 
 Then  $\beta_{X}\hskip-2.5pt\restriction_Z=\xi$ is the finest locally convex topology, so $\beta_Y$ is of the same type on $L_{p}(Y)$. 
Hence $C_{p}(Y)$ is distinguished, by Theorem \ref{three}, equivalently, $Y$ is a $\Delta$-space, by Theorem  \ref{Theor:description}.
\end{proof}

If $C_p(X)$ is homeomorphic to a retract of $\R^\kappa$ for some cardinal $\kappa$, then $X$ is discrete \cite[Problem 500]{Tkachuk_book2}.
Nevertheless, there exists a continuous mapping from $\R^\omega$ onto $C_p[0,1]$ \cite[Problem 486]{Tkachuk_book2}.
Several open problems have been posed in the following direction: Suppose that a dense subspace of $C_p(X)$ is a "nice" (not necessarily linear)
continuous image of $\R^\kappa$,
for some cardinal $\kappa$; must $X$ be discrete? \cite[Section 4.2]{Tkachuk_book1}.
Lemma \ref{help} implies immediately  

\begin{corollary}\label{discrete}
Let a dense subspace of $C_p(X)$ be a continuous linear image of $\R^\kappa$,
for some cardinal $\kappa$. Then $X$ is discrete.
\end{corollary}

For simplicity, a topological space $X$ is called a \emph{$Q$-space} if each subset of $X$ is $F_{\sigma}$, or, equivalently,
each subset of $X$ is $G_{\delta}$ in $X$.

\begin{theorem}\label{linear2}
Let $X$ and $Y$ be normal spaces and assume that $Y$ is $l$-dominated by $X$.
 If $X$ is a $Q$-space, then $Y$ also is a $Q$-space.
\end{theorem}
\begin{proof} Normal $X$ is a $Q$-space if and only $X$ is \emph{strongly splittable}, i.e. for every $f\in\R^X$ there exists a sequence 
$S=\{f_n: n\in\omega\} \subset C_p(X)$ such that $f_n \rightarrow f$ in $\R^X$, by \cite[Problems 445, 447]{Tkachuk_book3}.
Let $T: C_{p}(X) \longrightarrow C_{p}(Y)$ be a continuous linear surjection.
Denote by $\widehat{T}:\mathbb{R}^{X}\longrightarrow \mathbb{R}^{Y}$ the extension of $T$ which is supplied by Lemma \ref{help}.
Take arbitrary $f\in\mathbb{R}^{Y}$. There exists $g\in\mathbb{R}^{X}$ with $\widehat{T}(g)=f$. Then there exists a sequence $B\subset C_{p}(X)$
converging to $g$ in $\R^X$. 
We define $A = T(B)$. It is easy to see that $A \subset C_p(Y)$ converges to $f$ in $\R^Y$.
\end{proof}

\begin{corollary}\label{Q-set}
Let $X$ and $Y$ be metrizable spaces (in particular, subsets of $\R$) and assume
that $Y$ is $l$-dominated by $X$.
If $X$ is a $Q$-set, then $Y$ also is a $Q$-set.
\end{corollary}

\begin{remark} Note that Theorems \ref{linear1} and \ref{linear2}, and Corollary \ref{Q-set} are valid
under a weaker assumption that a dense subspace of $C_p(Y)$ is a continuous linear image of $C_p(X)$.
\end{remark}

A space $X = \bigcup_{n\in\omega} X_n$ is called \emph{$\sigma$-scattered} (\emph{$\sigma$-discrete}) if every
$X_n$ is scattered (discrete, respectively).

\begin{proposition}\label{linear3}
Assume that $Y$ is $l$-dominated by $X$.
If $X$ is $\sigma$-scattered ($\sigma$-discrete), then $Y$ also is $\sigma$-scattered ($\sigma$-discrete, respectively). 
\end{proposition}
\begin{proof} 
Our argument is a modification of the proof of \cite[Theorem 3.4]{LLP} and is based 
on an analysis of the dual spaces (see also \cite[Proposition 2.1]{Kawamura}).  
Recall that for a Tychonoff space $X$, $L_{p}(X)$ denotes the dual space, that is, the space of all continuous 
linear functionals on $C_{p}(X)$ endowed with the pointwise convergence topology. For each natural $n\in \N$
consider the subspace $A_n(X)$ of $L_{p}(X)$ formed by all words of the reduced length precisely $n$.
It is known that $A_n(X)$ is homeomorphic to a subspace of the Tychonoff product $(\R^*)^n\times X^n$, where $\R^*= \R\setminus\{0\}$. 
Let $T: C_{p}(X) \longrightarrow C_{p}(Y)$ be a continuous linear surjection.
The adjoint mapping $T^*$ embeds $L_p(Y)$ into $L_p(X)$. Therefore, $Y$ can be represented as a countable union of
subspaces $Y_i, i \in \N$, such that each $Y_i$ is homeomorphic to a subspace of $(\R^*)^n\times X^n$ for some $n=n(i)$.

Consider the projection $p_i$ of each of the above pieces $Y_i \subset (\R^*)^n\times X^n$ to the second factor $X^n$.
The surjectivity of the linear mapping $T$ implies that $p_i: Y_i \longrightarrow X^n$ is a finite-to-one mapping.
Evidently, $X^n$ is scattered /discrete provided $X$ is. Since $p_i$ is continuous, for every isolated point $z \in X^n$ its finite fiber 
$p_i^{-1}(z)$ consists of points isolated in $Y_i$ and the claim follows.
\end{proof}

The following question remained open.
\begin{problem} \label{prob5}
Let $Y$ be $l$-dominated by a scattered space $X$. Must $Y$ be scattered?
\end{problem}

Below we answer Problem \ref{prob5} positively in several interesting particular cases with the help of the properties of $\Delta$-spaces.

\begin{proposition}\label{prop:completely}
Let $X$ and $Y$ be metrizable spaces and assume that $Y$ is $l$-dominated by $X$.
If $X$ is scattered, then $Y$ also is scattered.
\end{proposition}
\begin{proof} If $X$ is metrizable and scattered, then $X$ is a $\Delta$-space by \cite[Proposition 4.1]{KL}.
Hence by Theorem \ref{linear1} the space $Y$  is a $\Delta$-space.
From another hand, every metrizable and scattered space is completely metrizable, by \cite[Corollary 2.2]{Tkachuk-1}.
A metrizable space $Y$ is $l$-dominated by a completely mertizable space $X$, therefore $Y$ is completely metrizable by the main result of \cite{Pelant}.
Finally, $Y$ is a \v{C}ech-complete $\Delta$-space, and $Y$ is scattered applying \cite[Theorem 3.4]{KL}. 
\end{proof}

If $X$ and $Y$ both are compact spaces and there is a continuous mapping from $C_p(X)$ onto $C_p(Y)$, then $Y$ is Eberlein whenever $X$ is (see \cite[Theorem IV.1.7]{Arch5}),
and $Y$ is Corson whenever $X$ is (see \cite[Theorem IV.3.1]{Arch5}). 
 Our next statement is a combination of a few known results,
while we apply Theorem \ref{linear1} in order to obtain the scatteredness of a target space.

\begin{proposition}\label{prop:Eber} 
Assume that $Y$ is $l$-dominated by $X$.
\begin{enumerate}
\item[{\rm (1)}] If $X$ is an Eberlein compact, then $Y$ also is an Eberlein compact.
\item[{\rm (2)}] If $X$ is a scattered Eberlein compact, then $Y$ also is a scattered Eberlein compact.
\end{enumerate} 
\end{proposition}
\begin{proof} (1)
The space $C_p(X)$ contains a dense $\sigma$-compact subspace, by \cite[Theorem IV.1.7]{Arch5}, hence $C_p(Y)$ satisfies the same property and  
consequently, $C_{p}(Y)$ contains a compact subset $K$ which separates points of $Y$.
On the other hand, $Y$ is pseudocompact by the result of V. Uspenskii (see \cite{Arch4}).
By means of evaluation mapping we define a continuous injective mapping $\varphi: Y  \longrightarrow C_p(K)$.
Denote by $B=\varphi(Y)$. Then $B$ is a pseudocompact subspace of $C_p(K)$. Applying \cite[Theorem IV.5.5]{Arch5} we get that $B$ is an Eberlein compact.
 We showed that the pseudocompact space $Y$ is mapped by a continuous injective mapping $\varphi$ onto the Eberlein compact $B$.
However, the mapping $\varphi$ must be a homeomorphism, by \cite[Theorem IV.5.11]{Arch5} and the result follows.

(2) Every scattered Eberlein compact is a $\Delta$-space, by \cite[Theorem 3.7]{KL}.
Hence by above Theorem \ref{linear1}, $Y$ is a $\Delta$-space. We conclude that $Y$ is scattered, by \cite[Theorem 3.4]{KL}.
\end{proof}

\begin{remark}\label{Rez}
Proposition \ref{prop:Eber}(1) is not valid for Corson compacts.
E. Reznichenko showed that there exists a compact space $X$ with
the following properties (see \cite{Arch5}, \cite[Problem 222]{Tkachuk_book3}):
\begin{enumerate}
\item[{\rm (i)}] $C_p(X)$ is a $K$-analytic space, i.e. $X$ is a Talagrand (hence, Corson) compact;
\item[{\rm (ii)}] there is $x \in X$  such that $Y= X \setminus \{x\}$ is pseudocompact and $X$ is the Stone-\v{C}ech compactification of $Y$.
\end{enumerate}
Evidently, the restriction mapping projects continuously $C_p(X)$ onto $C_p(Y)$.
\end{remark}

\begin{remark}\label{Mar}
The assumption of linearity of continuous surjection between function spaces, even for compact spaces $X$ and $Y$,
 cannot be dropped in the main Theorem \ref{linear1} and its corollaries above.
Let $Y$ be any non-scattered metrizable compact, for instance $Y = [0,1]$. Denote by $S$ the convergent sequence.
Using the argument presented in \cite[Proposition 5.4]{Krupski}, one can construct a continuous surjective mapping from $C_{p}(S)$ onto $C_{p}(Y)$ (see also \cite[Remark 3.4]{Kawamura}). 
\end{remark}

\begin{problem}\label{uniform}
 Assume that $Y$ is $u$-dominated by $X$.
 Is it true that $Y$ is a $\Delta$-space provided $X$ is a $\Delta$-space?
\end{problem}

Theorem \ref{linear1} may suggest also the following questions.
Below $C_{k}(X)$ stands for the space of all real-valued continuous functions on a Tychonoff space $X$
 endowed with the compact-open topology.
\begin{problem}\label{prob6}
 Assume that $X$ and $Y$ are Tychonoff spaces and there exists a continuous linear surjection from $C_{k}(X)$ onto $C_{k}(Y)$.
 Is it true that $Y$ is a $\Delta$-space provided $X$ is a $\Delta$-space?
\end{problem}

In case when the answer to Problem \ref{prob6} is negative one can pose the following

\begin{problem}\label{prob7}
Find scattered compact spaces $X$ and $Y$ such that $X\in\Delta$ but $Y\notin\Delta$ 
and there exists a continuous linear surjection from the Banach space $C(X)$ onto the Banach space $C(Y)$.
\end{problem}
Surely, for such $X$ and $Y$ a continuous linear surjection from $C_{p}(X)$ onto $C_{p}(Y)$ does not exist, by Theorem \ref{linear1}. 
Notice also that $X$ in Problem \ref{prob7} cannot be an Eberlein compact, since otherwise $Y$ would be a scattered Eberlein compact, hence $Y$ would be in the class $\Delta$.

\begin{proposition}\label{Cech-complete}
Let $X$ be a \v{C}ech-complete Lindel\"of space. Then
the following assertions are equivalent. 
\begin{enumerate}
\item[{\rm (1)}] $X$ is scattered.
\item[{\rm (2)}] $X$ is $\sigma$-scattered.
\item[{\rm (3)}] $C_{p}(X)$ is a Fr\'echet-Urysohn space.
\end{enumerate}
\end{proposition}
\begin{proof} The implication (2) $\longrightarrow$ (1) follows from the well-known fact that
every closed subspace of $X$ satisfies Baire category theorem.
The equivalence (1) $\longleftrightarrow$ (3) has been proved already in \cite[Corollary 2.12]{ga-ka}.
\end{proof}

\begin{corollary}
Let $X$ be a \v{C}ech-complete Lindel\"of space. If $X \in \Delta$, then $C_{p}(X)$ is a Fr\'echet-Urysohn space.
\end{corollary}

Let us call a lcs $E$ \emph{hereditarily distinguished} if every closed linear subspace of $E$ is distinguished.
It is known that even a Fr\'echet distinguished lcs can contain a closed non-distinguished subspace.
The only hereditarily distinguished $C_p$-spaces we are aware of are the products of reals $\R^{\kappa}$.
Note that if $\varphi$ is a continuous mapping from a compact $\Delta$-space $X$ onto $Y$, then
the adjoint mapping $\varphi^*$ identifies $C_p(Y)$ with a closed linear subspace of $C_p(X)$,
and this closed copy of $C_p(Y)$ is again distinguished. Our last problem is inspired by this observation.

\begin{problem}\label{prob8}
Does there exist an infinite compact space $X$ such that $C_p(X)$ is hereditarily distinguished?
More specifically, let $X$ be the one-point compactification of an infinite discrete space. Is $C_p(X)$ hereditarily distinguished?
\end{problem}

\section{Compact $\Delta$-spaces and PFA}\label{PFA}
A topological space $X$ has countable tightness if for each $A \subset X$ and for each $x \in \cl(A)$, there is a countable $B \subset A$ such that $x \in \cl(B)$.
A topological space $X$ is a sequential space if $A \subset X$ is a sequentially closed set implies that $A$ is closed.
 The set $A$ is a sequentially closed set if a countable sequence $(x_n) \in A$ converges to $x \in X$ implies that $x \in A$.
Every sequential space is countably tight. 
A topological space $X$ is called \emph{$\omega$-bounded} if the closure of every countable subset of $X$ is compact.
$X$ is called \emph{pseudocompact} if every continuous function defined on $X$ is bounded.
Evidently, every $\omega$-bounded space is countably compact, and every countably compact space is pseudocompact.
The space of countable ordinals $[0,\omega_{1})$ is an example of an $\omega$-bounded space which is not compact.
 A continuous mapping $f: X \longrightarrow Y$ is called perfect if  
it is closed and $f^{-1}(y)$ is compact for each $y\in Y$.

\begin{theorem}\label{theor:omega}
Every $\omega$-bounded $\Delta$-space is compact. 
\end{theorem} 
\begin{proof} Assume that $X$ is a counter-example to the claim. Then, by a result of D. Burke and G. Gruenhage \cite[Lemma 1]{Gruenhage}, 
$X$ contains a subset $Z$ which is a perfect preimage of the ordinal space $[0,\omega_{1})$. We conclude that a $\Delta$-space $Z$ can be mapped by a continuous closed mapping onto 
$[0,\omega_{1})$. By our Proposition \ref{prop:closedmap} this would mean that $[0,\omega_{1}) \in \Delta$,  however, the opposite is true \cite{KL}.
Obtaining contradiction finishes the proof.
\end{proof}

It has been shown in \cite{LS} that every compact $\Delta$-space has countable tightness.
Essentially the same argument as in Theorem \ref{theor:omega} provides a very short proof of this assertion.

\begin{theorem}\label{tightness}\cite{LS}
Every compact $\Delta$-space has countable tightness.
\end{theorem}
\begin{proof} A compact space has countable tightness if and only if it does not contain a perfect preimage of $[0,\omega_{1})$ (see \cite{Balogh}).
We argue again that every $\Delta$-space satisfies this property.
\end{proof}

A very natural question arises whether Theorems \ref{theor:omega} and \ref{tightness} can be generalized for countably compact spaces.
A positive answer follows from the Proper Forcing Axiom (PFA), due to the celebrated results of Z. Balogh \cite{Balogh} (see also \cite{BDFN}).

\begin{theorem}\label{theor:PFA} 
{\rm (PFA)}
\begin{enumerate}
\item[{\rm (1)}] Every countably compact $\Delta$-space is compact.
\item[{\rm (2)}] Every countably compact $\Delta$-space has countable tightness.
\item[{\rm (3)}] Every countably compact $\Delta$-space (hence, every compact $\Delta$-space) is sequential.  
\end{enumerate}
\end{theorem}

\begin{problem}\label{prob9}
Is it possible to obtain the results of Theorem \ref{theor:PFA} in ZFC alone?
\end{problem}

Note that all known examples of compact $\Delta$-spaces are $\sigma$-discrete. However, we don't know
if it is always the case.

\begin{problem}\label{prob3}
Let $X$ be a compact $\Delta$-space. Is $X$ a $\sigma$-discrete space?
\end{problem}

A closely related question to the last problem is the following one: When a $\Delta$-space is scattered?
 As we have been mentioned earlier every \v Cech-complete $\Delta$-space
 (in particular, every compact $\Delta$-space) is scattered \cite{KL}.

\begin{example}\label{ex:Baire} There exists a Baire countable space which is not scattered.
Fix in the real line $\R$ a countable dense subset $B$ consisting of irrationals. Let $X$ be the union of the rationals $\Q$ with $B$.
Equip $X$ with the topology inherited from the Michael line $M$. Then $X$ is a countable space containing a copy of $\Q$, therefore $X$ is a non-scattered $\Delta$-space.
$X$ is Baire since it contains a dense discrete subspace.
\end{example}

Despite of Problem \ref{prob9} the following result does not require extra set-theoretic assumptions.

\begin{theorem}\label{th:count_compact} Every countably compact $\Delta$-space is scattered.
\end{theorem}
\begin{proof} On the contrary, assume that a countably compact space $X$ is not scattered.
Every countably compact space $X$ is pseudocompact, therefore there exists a closed subset $K \subset X$ and a continuous surjective mapping $\varphi$
from $K$ onto the closed interval $[0,1]$, by \cite[Proposition 5.5]{LT}. Every closed subset $F$ of $K$ is a countably compact space, its continuous image
$\varphi(F)$ is a countably compact subset of $[0,1]$, therefore $\varphi(F)$ is compact. We conclude that $\varphi\,\res\,{K}$ is a closed continuous mapping from $K$ onto $[0,1]$.
This evidently contradicts Theorem \ref{th:F-sigma_map}, since $[0,1]\notin\Delta$. 
\end{proof}

The proof above fails if we assume only that $X$ is a pseudocompact (and non-normal) space, in view of Proposition \ref{prop:locomp}.

\begin{problem}\label{prob4} Prove that every pseudocompact $\Delta$-space is scattered, thereby prove that a Tychonoff space $Y$ is scattered provided it is 
$l$-dominated by a compact $\Delta$-space $X$.
\end{problem}
\textbf{Acknowledgments.} The authors thank Witold  Marciszewski for the useful information about Corson and Eberlein compact spaces.

The authors acknowledge and thank Vladimir Tkachuk for the most stimulating letters \cite{Tkachuk}.


\end{document}